%%%%%%%%%%%%%%%%%%%%%%%%%%%%%%%%%%%%%%%%%%%%%%%%%%%%%%%%%%%%%%%%%%%%%%%%%%%
%% Lohkamp, Joachim
%% 
%% Negatively Ricci curved manifolds
%% 
%% In this paper we announce the following result: ``Every manifold of 
%%   dimension $\ge3$ admits a complete negatively Ricci curved metric.'' 
%%   Furthermore we describe some sharper results and sketch proofs.
%% 
%% publ:  Bull. Amer. Math. Soc. (N.S.) 27(1992) no. 2
%% pp:    288-291
%% type:  Research Announcement        markup: amstex    file size: 14K
%% 
%% copyright: American Math. Society copyright; see end of article
%% 
%% Include files necessary for this article: bull-ppt.tex
%% 
%%%%%%%%%%%%%%%%%%%%%%%%%%%%%%%%%%%%%%%%%%%%%%%%%%%%%%%%%%%%%%%%%%%%%%%%%%%
\input amstex %version 1.1
\documentstyle{amsppt}
\input bull-ppt
\keyedby{bull325/pxp}
\define\Diff{\operatorname{Diff}}
\define\Isom{\operatorname{Isom}}
\define\id{\operatorname{id}}

\topmatter
\cvol{27}
\cvolyear{1992}
\cmonth{October}
\cyear{1992}
\cvolno{2}
\cpgs{288-291}
\title Negatively Ricci curved manifolds\endtitle
\author Joachim Lohkamp\endauthor
\shortauthor{Joachim Lohkamp}
\shorttitle{Negatively Ricci curved manifolds}
\address Ruhr-Universit\"at Bochum, Mathematisches 
Institut, D-4630 Bochum, 
Germany\endaddress
\date November 6, 1991\enddate
%\keywords{}
\abstract In this paper we announce the following result: 
``Every manifold of
dimension $\ge3$ admits a complete negatively  Ricci 
curved metric.''
Furthermore we describe some sharper results and sketch 
proofs.\endabstract
\endtopmatter

\document
\heading I. Introduction\endheading
There are three well-known curvatures: sectional, Ricci, 
and scalar curvature.
While the existence of complete negatively sectional 
curved metrics leads to
many topological implications for the underlying manifold, 
each manifold admits
a complete metric of constant negative scalar curvature.

Since Ricci curvature takes position between these two 
curvatures, it is
reasonable to look for obstructions as well as for 
existence results for (complete)
negatively Ricci curved metrics.

This led to some long-standing conjectures in Riemannian 
geometry and appears
in lists of problems compiled by Yau \cite{Y}, problem 24, 
Kazdan \cite{K},
problem 9, Bourguignon \cite{Bg}, question 4.11, and others.

A first notable result was proved by Gao and Yau (cf. 
\cite{G, GY}). They
started from Thurston's Hyberbolic Dehn Surgery \cite{T}: 
Every compact
three-manifold $M$ can be obtained from $S^3$ by Dehn 
Surgery, along some link
$L_M\subset S^3$, which complement $S^3\backslash L_M$, 
admits a complete
hyperbolic metric with finite volume. Gao and Yau managed 
to perform this Dehn
Surgery such that the Ricci curvature $r(g)$ remains 
negative near $L_M$ and is
just the hyperbolic metric outside, i.e., one gets a 
metric with $r(g)<0$ on
each compact three-manifold. Finally they extended this 
result to
three-manifolds of finite type (with complete metric). 

Subsequently, there was a paper, written by Brooks 
\cite{Br}, that addressed
general existence theorems. He develops a technique of 
smoothing hyperbolic
orbifold singularities of order $K\ge 12$, which could be 
locally realized as
the quotient of hyperbolic space by an element of order 
$K$ fixing a
codimension 2 manifold $N$ (which is hyperbolic if $\dim 
N\ge 2)$. 

On the other hand, Thurston theory yields such higher 
singular hyperbolic
orbifold metrics on each compact three-manifold and Brooks 
combines these ideas
to get: Every compact three-manifold admits a metric $g$ 
with $-a<r(g)<-b$ for
constants $a>b>0$, which are independent of the chosen 
manifold.

In contrast to that we will develop a completely different 
method of
attack. 
Starting from an almost arbitrary metric, we obtain our 
result by local
deformations, which will be sketched in part III of this 
paper.

\heading II. New results\endheading
\thm{Theorem 1} There are constants $a(n)>b(n)>0$ 
depending only on the
dimension $n\ge 3$, such that each manifold $M^n$ admits a 
complete metric $g$
with 
$$
-a(n)<r(g)<-b(n)\.
$$
\ethm

\noindent Notice that even for $n=3$ we give a proof that 
does not use Thurston theory.

Bochner proved a classical result: The isometry group of a 
compact
manifold with $r(g)<0$ is finite. We will show that this 
``geometric
restriction'' is sharp.

\thm{Theorem 2} If $M^n$, $n\ge 3$, is a compact manifold 
and $G$ is a subgroup
of $\Diff(M)$, then
$$
G=\Isom(M,g)\quad\text{for some metric }g\text{ with 
}r(g)<0\Leftrightarrow
G\text{ is finite}\.
$$
\ethm

Next, there are unexpected density and cut-off properties 
of negatively
Ricci-curved metrics.

Gromov \cite{Gr} introduced the so-called Hausdorff 
distance $d_{\roman H}$ 
between two metric spaces $M_1$, $M_2$, which can be 
viewed as the minimal
distance or deviation between $M_1$ and $M_2$ for all 
possible isometric
embeddings in any metric space $M$. Of course, Riemannian 
manifolds can also be
considered as metric spaces, and we get for compact $M^n$, 
$n\ge 3$,

\thm{Theorem 3} $\{(M^n,g)|r(g)<0\}$ is dense in 
$\{(M^n,g)|g$ arbitrary
metric$\}$ with respect to $d_{\roman H}$.
\ethm

Our final result is the most flexible one. It is a useful 
tool for constructing
negatively curved metrics.

\thm{Theorem 4} Let $S\subset M^n$, $n\ge 3$, be a closed 
subset, $S\subset U$
an open neighborhood, and $g_0$ a metric on $S$ with 
$r(g_0)\le0$ \RM(resp.
$r(g_0)<0)$ on $U$. Then there is metric $g$ on $M$ with 
$r(g)\le 0$ \RM(resp.
$r(g)<0)$ on $M$, and $g\equiv g_0$ on $S$.
\ethm

\rem{Remark} Theorem 1 implies that each manifold admits a 
complete metric with
scalar curvature bounded by two negative constants. In 
this case, the Yamabe
equation can be solved easily, and we get a complete 
metric of constant
negative scalar curvature. This was proved before by Aubin 
\cite{A} and 
Bland and Kalka \cite{BK} in a different way.
\endrk

\heading III. Some ideas of the proof\endheading
The proof is by induction (with respect to dimension) and 
consists of two main
parts:
\roster
\item There is a metric $g_n^-$ on $\Bbb R^n$ with 
$r(g_n^-)<0$ on $B_1(0)$ and
$g_n^-\equiv g_{Eucl.}$ on $\Bbb R^n\backslash B_1(0)$.
\item Using this metric $g_n^-$ a Besicovitch type 
covering argument yields the
theorems for $n$-dimensional manifolds.
\endroster

In part (1) there are two cases where the proofs are 
completely different:
$n=3$ and $n\ge 4$.

The case $n=3$ starts from a metric $g$ on $\Bbb 
R^3\bigsharp S^1\times S^2$,
which is isometric to $(\Bbb R^3\backslash 
B_1(0),g_{Eucl.})$ outside a set
$U$, with $r(g)<0$ on $U$. By some warped product trick, 
which also allows to
extend Brooks results to general hyperbolic singularities 
(cf. part I.), we obtain
$g_3^-$. In this short note we cannot go into further 
details. Instead (for the
moment) we will use Gao and Yau's result: $S^3$ admits a 
metric with $r(g)$. So
we easily obtain $g_3^-$ by taking connected sums $S^3 
\;\#\;\Bbb R^3=\Bbb R^3$.

The proof of case $n\ge4$ is more complicated: The main 
step consists in the
construction of a metric $g(n)$ on $S^1\times 
B_2(0)\subset S^1\times\Bbb
R^{n-1}$ with $r(g(n))<0$ on $S^1\times B_1(0)$ and 
$g(n)=g_{S^1}+g_{Eucl.}$
outside. To get this, one embeds $\overline{S^1}\times 
S^1\times
B_5^{n-2}(0)\subset\overline{S^1}\times S^1\times\Bbb 
R^{n-2}$ into $S^1\times
B_{1/2}(0)$ in some special way and uses Theorem 4 in 
dimension $n-1$ for
$\overline{S^1}\times  B_5^{n-2}(0)$ and $S^1\times 
B_5^{n-2}(0)$ to construct two metrics
the combination of which yields (after some minor 
manipulations) $g(n)$.
Finally, one embeds $S^1\times B_3(0)$ into $\Bbb R^n$ as 
tubular neighborhood
of a large circle. Some not-too-hard deformation argument 
yields $g_n^-$. 

To give an idea of part (2), we sketch the existence proof 
for negatively Ricci
curved metrics in a very special case: the $n$-dimensional 
flat torus
$T^n=S^1\times\cdots\times S^1$, each factor with length 
$2\pi\cdot 100$. For each
$\varrho\in\ ]0,1[$ there is a discrete subset 
$A_\varrho\subset T^n$ with
\roster
\item "(i)" \<$d(a,b)>5\cdot\varrho$ for $a\ne b\in 
A_\varrho;$
\item "(ii)" \<$T^n=\bigcup_{a\in 
A_\varrho}\overline{B_{5\cdot\varrho}(a)}$;
\item "(iii)" \<$\bigsharp\{a\in A_\varrho|z\in 
B_{10\cdot\varrho}(a)\}\le
c(n)$, $c(n)$ independent of $z\in T^n$ and $\varrho\in 
]0,1[$.
\endroster

We define $g_{A_\varrho}\coloneq \varrho^2\cdot 
f_{a,\varrho}^*(g_n^-)$ on
$B_{2\cdot\varrho}(a)$, $a\in A_\varrho$, and 
$g_{A_\varrho}\equiv g_{T^n}$
elsewhere with $f_{a,\varrho}\:B_{2\cdot\varrho}(a)\to\Bbb 
R^n$,
$f_{a,\varrho}\equiv(1/\varrho)\cdot I_a\circ\exp_a^{-1}$, 
where
$I_a\:T_aT^n\to \Bbb R^n$ is some linear isometry 
$(g_{A_\varrho}$ depends on
the choice of $I_a$; for Theorem 2 one chooses 
$I_{f(a)}^{-1}\circ I_a\equiv
Df_a$ for each $f\in G\subset\Diff(M))$. We are now ready 
to consider the
following metric:
$$
g(A_\varrho,d,s)\coloneq\prod_{a\in 
A_\varrho}\exp(2F^\varrho_{d,s}\cdot
h_\varrho(10\cdot\varrho-d(a,\id_{T^n})))\cdot g_{A_\varrho}
$$
where $F^\varrho_{d,s}\equiv 
s\cdot\exp(-d\cdot\varrho/\id_{\Bbb R})$, $h_\varrho\equiv
h((1/\varrho)\cdot\id_{\Bbb R})$ with $h\in C^\infty(\Bbb 
R,[0,1])$,
$h\equiv 0$ on $\Bbb R^{\le 1/2}$, $h\equiv1$ on $\Bbb 
R^{\ge 3/4}$, and
$d(\cdot,\cdot)$ is the usual distance on $T^n$.

Then there are $d_0,s_0>0$ independent of $\varrho$ such 
that
$r(g(A_\varrho,d,s))<0$ for each $d>d_0$, $s\in\ ]0,s_0[$. 
Of course, for the
flat torus $T^n$, $\varrho$ is a useless parameter, but 
think of the general
case: A Besicovitch type argument gives similar coverings 
for each manifold
$M^n$ (with $c(n)$ independent of $M^n)$ if a suitable 
start metric is chosen.
If $\varrho$ is chosen very small, the background metric 
on the given manifold
appears quite ``flat relative to $B_\varrho(a)$\<''. 
Again, there are
$\varrho,d,s$ such that an analogously defined metric 
$g(A_\varrho,d,s)$ is
negatively Ricci curved, but $\varrho,d$, and $s$ are no 
longer independent.

\heading Acknowledgment\endheading
This paper is an abstract of the author's doctoral thesis. 
The author thanks
Professor J\"urgen Jost for his friendly support.


\Refs
\ref\key A\by T. Aubin\paper M\'etriques riemanniennes et 
courbure
\jour J. Differential Geom.\vol 4\yr 1970\pages 
383--424\endref
\ref\key BK\by J. Bland and M. Kalka\paper Negative scalar 
curvature metrics on
noncompact manifolds\jour Trans. Amer. Math. Soc.\vol 
316\yr 1989\pages 433--446\endref
\ref\key Bg\by J. P. Bourguignon\paper Ricci curvature and 
Einstein metrics,
{\rm Global Differential Geometry}\inbook Lecture Notes in 
Math., vol. 838
\publ Springer, New York, 1981\pages 42--63\endref 
\ref\key Br1\by R. Brooks\paper A construction of metrics 
of negative Ricci
curvature\jour J. Differential Geom.\vol 29\yr 1989\pages 
85--94\endref
\ref\key G\by L. Z Gao\paper The construction of 
negatively Ricci curved
manifolds\jour Math. Ann.\vol 271 {\rm(1985)}\pages 
185--208\endref
\ref\key GY\by L. Z. Gao and S. T. Yau\paper The existence 
of negatively Ricci
curved on three manifolds\jour Invent. Math.\vol 85\yr 
1986\pages
637--652\endref
\ref\key Gr\by M. Gromov\book Structures m\'etriques pour 
les vari\'et\'es
riemanniennes\publ Editions CEDIC\publaddr Paris\yr 
1981\endref
\ref\key K\by J. L. Kazdan\paper Prescribing the curvature 
of a Riemannian
manifold\inbook CBMS Regional Conf. Ser. in Math., vol. 
57\publ Conf. Board Math.
Sci.\publaddr Washington, DC\yr 1985\endref
\ref\key T\by W. Thurston\paper The geometry and topology 
of three
manifolds\inbook Princeton Lecture Notes, vol. 57\publaddr 
Princeton, NJ\yr 1980\endref
\ref\key Y\by S. T. Yau\paper Seminar on differential 
geometry, problem
section\inbook Ann. of Math. Stud., vol. 102\publ 
Princeton Univ. Press\publaddr
Princeton, NJ\yr 1982\endref
\endRefs
\enddocument